% An algorithm to determine the Heegaard genus of simple 3-manifolds

\input amssym.def
\input amssym
\input epsf
\magnification=1100
\baselineskip = 0.23truein
\lineskiplimit = 0.01truein
\lineskip = 0.01truein
\vsize = 8.5truein
\voffset = 0.2truein
\parskip = 0.10truein
\parindent = 0.3truein
\settabs 12 \columns
\hsize = 5.8truein
\hoffset = 0.2truein

\setbox\strutbox=\hbox{%
\vrule height .708\baselineskip
depth .292\baselineskip
width 0pt}
\font\caps=cmcsc10

\font\bigtenrm=cmr10 at 14pt

\def\sqr#1#2{{\vcenter{\vbox{\hrule height.#2pt
\hbox{\vrule width.#2pt height#1pt \kern#1pt
\vrule width.#2pt}
\hrule height.#2pt}}}}
\def\square{\mathchoice\sqr46\sqr46\sqr{3.1}6\sqr{2.3}4}

\centerline{\bigtenrm AN ALGORITHM TO DETERMINE}
\centerline{\bigtenrm THE HEEGAARD GENUS OF SIMPLE 3-MANIFOLDS}
\centerline{\bigtenrm WITH NON-EMPTY BOUNDARY}
\tenrm
\vskip 14pt
\centerline{MARC LACKENBY}
\vskip 18pt
\centerline{\caps 1. Introduction}
\vskip 6pt

The Heegaard genus of a compact orientable 3-manifold is an important invariant.
The aim of this paper is to demonstrate that it is algorithmically
computable, at least when the 3-manifold is simple and has non-empty boundary.
Recall that a compact orientable 3-manifold is {\sl simple} if it
is irreducible and any properly embedded disc, incompressible annulus
or incompressible torus is boundary parallel.

\noindent {\bf Theorem 1.1.} {\sl Let $M$ be a compact connected orientable simple
3-manifold with non-empty boundary. Then there is
an algorithm to determine the Heegaard genus of $M$. Moreover,
for any given positive integer $n$, there is an algorithm
to find all Heegaard surfaces for $M$ with genus at most $n$
(up to ambient isotopy).}

This theorem can be applied to determine the tunnel number of hyperbolic
links. Recall that a {\sl tunnel system} for a link $L$ in $S^3$ is a collection
of disjoint embedded arcs $t$ with $t \cap L = \partial t$, such that the
exterior of $L \cup t$ is a handlebody. The {\sl tunnel number} of $L$
is the minimal number of arcs in a tunnel system. Two tunnel systems
$t_1$ and $t_2$ for $L$ are {\sl slide-equivalent} if there is an
isotopy of $S^3$ keeping $L$ fixed throughout, taking $N(L \cup t_1)$
to $N(L \cup t_2)$.

\noindent {\bf Corollary 1.2.} {\sl Let $L$ be a hyperbolic link in the 3-sphere.
Then there is an algorithm to determine the tunnel number of $L$.
Moreover, for any given positive integer $n$, there is an algorithm
to find all tunnel systems for $L$ with at most $n$ arcs
(up to slide-equivalence).}

The input to the algorithms in Theorem 1.1 is a triangulation of $M$.
In Corollary 1.2, one may supply a diagram of the link or a triangulation of
its exterior. Thus, the hyperbolic structure does not need to be given in advance.
The second algorithm provided by Theorem 1.1 creates a finite list of
Heegaard surfaces for $M$. More specifically, it provides an explicit
subdivision of the triangulation of $M$ and explicit subcomplexes
which are the required Heegaard surfaces. Note, however,
there is no guarantee that the surfaces in this
list are pairwise non-isotopic. This is because there is currently
no known algorithm for determining whether two Heegaard surfaces for
a 3-manifold are ambient isotopic.

Most of the key ideas behind this paper are due to Rubinstein.
He proved that, given any triangulation of a compact orientable
3-manifold, any strongly irreducible Heegaard surface may
be ambient isotoped into almost normal form [12]. Using the computable
nature of normal surface theory, he explained how one might 
use this to compute the Heegaard genus of the manifold.
However, the possible presence of normal tori creates formidable
technical obstacles to this approach. Jaco and Rubinstein [2]
have developed a theory of `1-efficient' and `layered' triangulations to
try to overcome these difficulties, but this appears to be highly
technical, and the results are not fully published. An alternative
approach to Heegaard surfaces has been developed by Li ([7],[8]),
starting with almost normal surfaces, but then using branched
surfaces. Using this theory, he has solved some important 
longstanding problems. One of his theorems is as follows.

\noindent {\bf Theorem 1.3.} [7] {\sl Any closed orientable irreducible
atoroidal 3-manifold has only finitely many Heegaard splittings
of a given genus, up to ambient isotopy.}

However, his proof is non-constructive, and so there appears
to be no immediate way of finding all these Heegaard surfaces using his techniques.
Our methods provide a proof of this result,
but where the 3-manifold is compact, connected, orientable and simple and has non-empty boundary.

The algorithms given in this paper follow Rubinstein's
original outline in many respects. Like Jaco and Rubinstein's
approach, the key is to use triangulations with very restricted normal
tori. But unlike their theory of 1-efficiency, the technique
here is to use angle structures. We introduce `partially
flat angled ideal triangulations', which have the key property
that they contain no normal tori other than those that
are normally parallel to a boundary component. We will
show that any compact connected orientable simple 3-manifold with
non-empty boundary (other than a 3-ball) has one of these ideal triangulations
and that there is an algorithm to construct it.

In a paper such as this, it is particularly important to be clear
about which parts are new and which are due to other mathematicians.
The material in Section 2, where partially flat angled ideal triangulations
are introduced, is new. However, similar notions have been used by other
authors for other purposes (for example [11]). Theorem 2.2,
which asserts that any compact connected orientable simple 3-manifold
with non-empty boundary (other than a 3-ball) has a partially flat 
angled ideal triangulation and that this may be
algorithmically constructed, is new. Section 3 contains
mostly expository material relating to generalised Heegaard
splittings. However, there are a number of important facts
in this section which appear in print for the first time.
These include Proposition 3.1, which states that, when one amalgamates a generalised
Heegaard splitting, the resulting Heegaard splitting is independent
of the choices that have been made. Additionally, we show that
if the generalised Heegaard splitting is given, say, as a 
subcomplex of a triangulation of the 3-manifold, then the
resulting Heegaard splitting is algorithmically constructible.
In Section 4, we state that a generalised Heegaard splitting
can be placed in normal and almost normal form, provided its
even surfaces are incompressible and have no 2-sphere components and its odd surfaces are
strongly irreducible. This is a mild
generalisation of a well-known result of Rubinstein [12] and
Stocking [15], and has essentially the same proof. We then describe the
computational aspects of normal and almost normal surfaces in
partially flat angled ideal triangulations. This is largely
routine. In the final section, we draw these many threads together
and describe the algorithms of Theorem 1.1.

\vskip 18pt
\centerline{\caps 2. Partially flat angled ideal triangulations}
\vskip 6pt

Angled ideal triangulations were first studied by Casson (unpublished),
and then developed by the author in [6]. They are just an ideal
triangulation, with an assignment of a real number to each
edge of each ideal tetrahedron, satisfying some simple conditions.
A mild generalisation of this concept, which we call a
partially flat angled ideal triangulation, is a
key ingredient of this paper.

An {\sl ideal tetrahedron} is a tetrahedron with its
vertices removed. An {\sl ideal triangulation} of a 
3-manifold $M$ is an expression of the interior of $M$ as a union of
ideal tetrahedra with their faces glued homeomorphically
in pairs. An {\sl angled ideal triangulation}
is an ideal triangulation, with a real number in the
range $(0, \pi)$ assigned to each edge of each ideal
tetrahedron, known as the {\sl interior angle} of the edge,
satisfying the following conditions:
\item{(i)} the angles at each ideal vertex of each
ideal tetrahedron sum to $\pi$;
\item{(ii)} the angles around each edge sum to
$2 \pi$.

In partially flat angled ideal triangulations, we 
allow some ideal tetrahedra to be {\sl flat}. This means
that the ideal tetrahedron is as shown in Figure 1. If two faces
of a flat ideal tetrahedron share an edge with interior
angle $\pi$, we term them {\sl coherent}. Thus,
the four faces are partitioned into two coherent
pairs.

\vskip 18pt
\centerline{
\epsfxsize=2.2in
\epsfbox{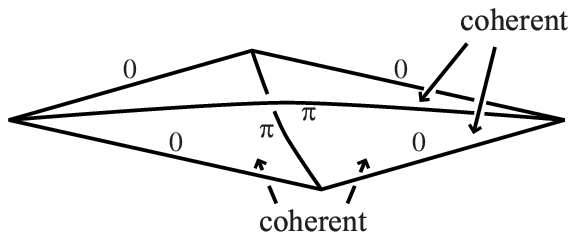}
}
\vskip 6pt
\centerline{Figure 1: a flat ideal tetrahedron}

A {\sl layered polygon} is a collection of flat ideal tetrahedra
glued together in a certain way to form a 3-manifold.
It is determined by the following data: an ideal polygon
with an initial ideal triangulation, together with a
sequence of elementary moves applied to this triangulation,
subject to the condition that every edge of the
initial triangulation that is not in the boundary of the
ideal polygon has a move performed upon it
at some stage. Recall that an {\sl elementary move}
on an ideal triangulation of a surface removes an
edge adjacent to two distinct ideal triangles, forming
an ideal square, and then inserts the other diagonal of
this square as a new edge. (See Figure 2.)

\vskip 18pt
\centerline{
\epsfxsize=1.8in
\epsfbox{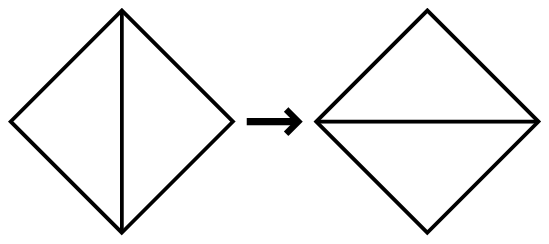}
}
\vskip 6pt
\centerline{Figure 2: an elementary move}

Starting with this data, we build
the layered polygon. We start with the initial ideal triangulation of
the ideal polygon, which will be the base of the layered polygon. 
The first move acts upon a pair of adjacent faces.
Attach onto them a flat ideal tetrahedron along a
coherent pair of faces. The `top' of the resulting
object inherits the second ideal triangulation of
the ideal polygon. Repeat this for each move of the
sequence. The resulting 3-manifold is the {\sl layered polygon}.
(See Figure 3.) It is a 3-ball with a finite collection of
points in its boundary removed.
Its boundary is the union of two ideal polygons, which
are the initial and terminal ideal polygons in the
sequence of elementary moves.
The intersection of these is a collection of edges,
which we term its {\sl vertical boundary}.

\vskip 18pt
\centerline{
\epsfxsize=4in
\epsfbox{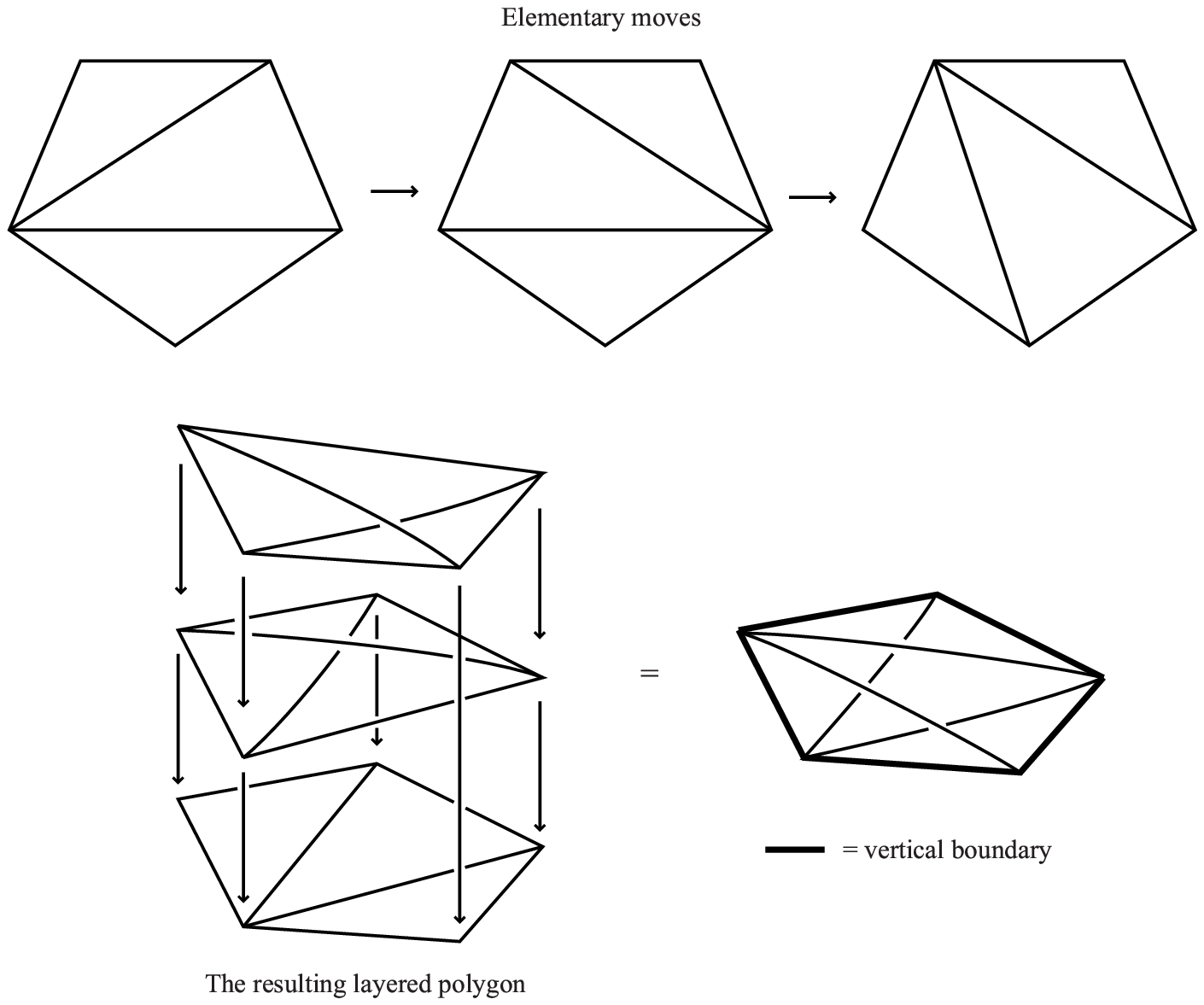}
}
\vskip 6pt
\centerline{Figure 3: construction of a layered polygon}

A {\sl partially flat angled ideal triangulation} of a 3-manifold
is an ideal triangulation, with a real number in the
range $[0, \pi]$ assigned to each edge of each ideal
tetrahedron, known as the {\sl interior angle} of the edge, satisfying the following conditions:
\item{(i)} the angles at each ideal vertex of each
ideal tetrahedron sum to at most $\pi$;
\item{(ii)} the angles around each edge sum to
$2 \pi$;
\item{(iii)} if the angles of an ideal tetrahedron are not
all strictly positive, then the ideal tetrahedron is flat;
\item{(iv)} the union of the flat ideal tetrahedra is a collection
of layered polygons, possibly with some edges in their vertical
boundary identified.

Note that, in (i), the angles at each ideal vertex are not
required to sum to precisely $\pi$, unlike the case of an
angled ideal triangulation. This is so that we can deal with 
3-manifolds having some boundary components with negative Euler characteristic.

Note also that we do not allow layered polygons to intersect each
other or themselves along anything other than vertical boundary edges.
They are not allow to touch at any point in the interior
of the `top' or `base' of a layered polygon. (See Figure 4.)

\vskip 18pt
\centerline{
\epsfxsize=2in
\epsfbox{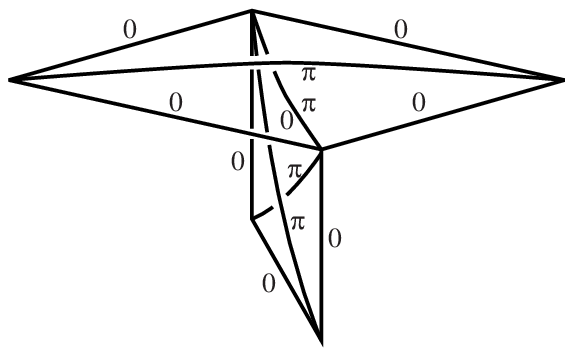}
}
\vskip 6pt
\centerline{Figure 4: a forbidden arrangement}

The usefulness of partially flat angled ideal triangulations
is that the normal and almost normal surfaces with non-negative
Euler characteristic that they contain are very constrained.
We briefly recall the relevant terminology.

A {\sl normal disc} in a tetrahedron or ideal tetrahedron is a properly embedded disc
that misses the vertices, that hits each
edge transversely in at most one point and that is not disjoint
from the edges. There are two
types of normal discs, {\sl triangles} and {\sl squares},
which are shown in Figure 5. A closed surface properly
embedded in $M$ is {\sl normal} if it intersects
each ideal tetrahedron in a collection of disjoint normal discs.

\vskip 18pt
\centerline{
\epsfxsize=3.3in
\epsfbox{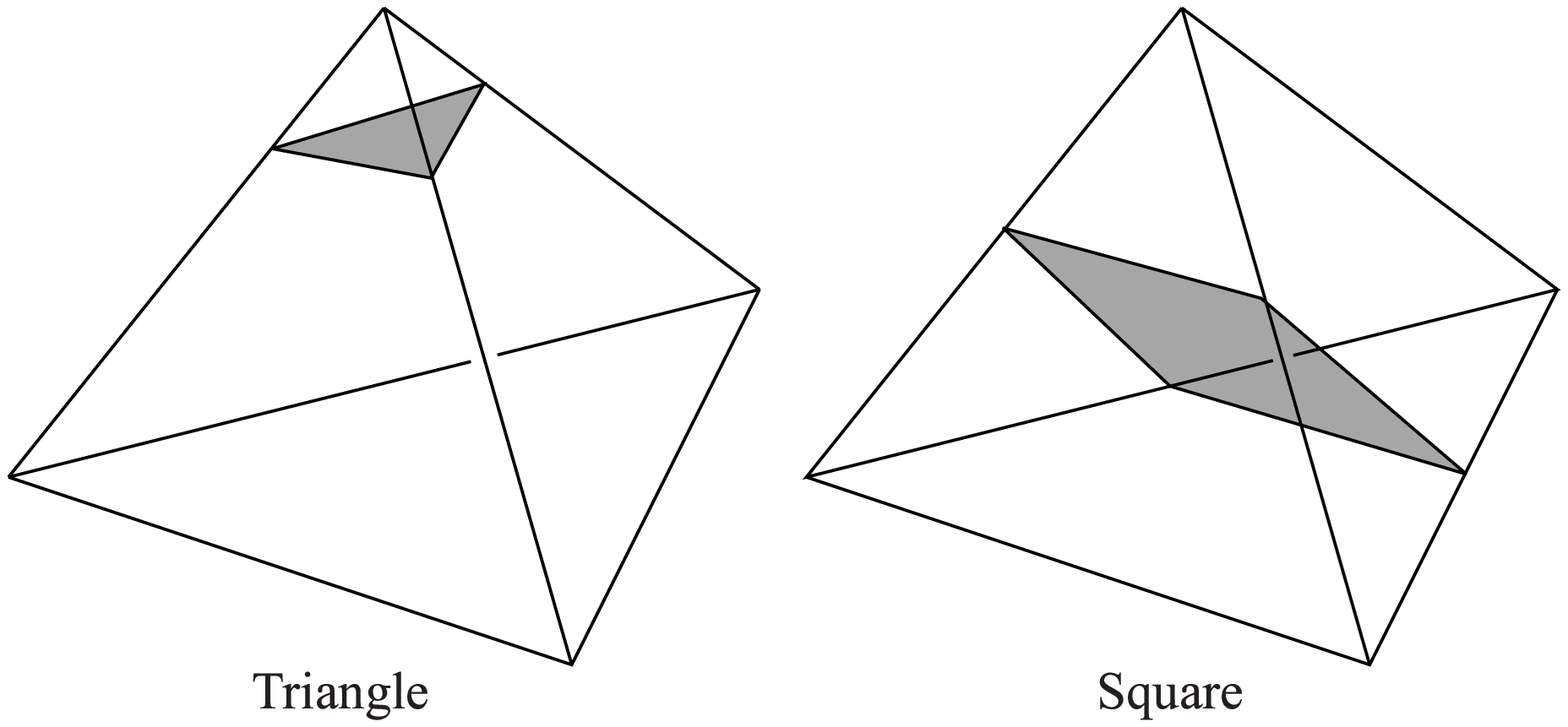}
}
\vskip 6pt
\centerline{Figure 5: normal discs}

An {\sl almost normal piece} in a tetrahedron or ideal tetrahedron is one of
two types: either an {\sl octagon}, as shown in Figure 6, or
a {\sl tubed piece}, which is two disjoint normal discs tubed together via a tube that runs parallel
to an edge. A closed surface properly embedded in
$M$ is {\sl almost normal} if it intersects each ideal
tetrahedron in a collection of normal discs, except
in precisely one ideal tetrahedron, where it is a collection
of normal discs and exactly one almost normal piece.

\vskip 18pt
\centerline{
\epsfxsize=3.3in
\epsfbox{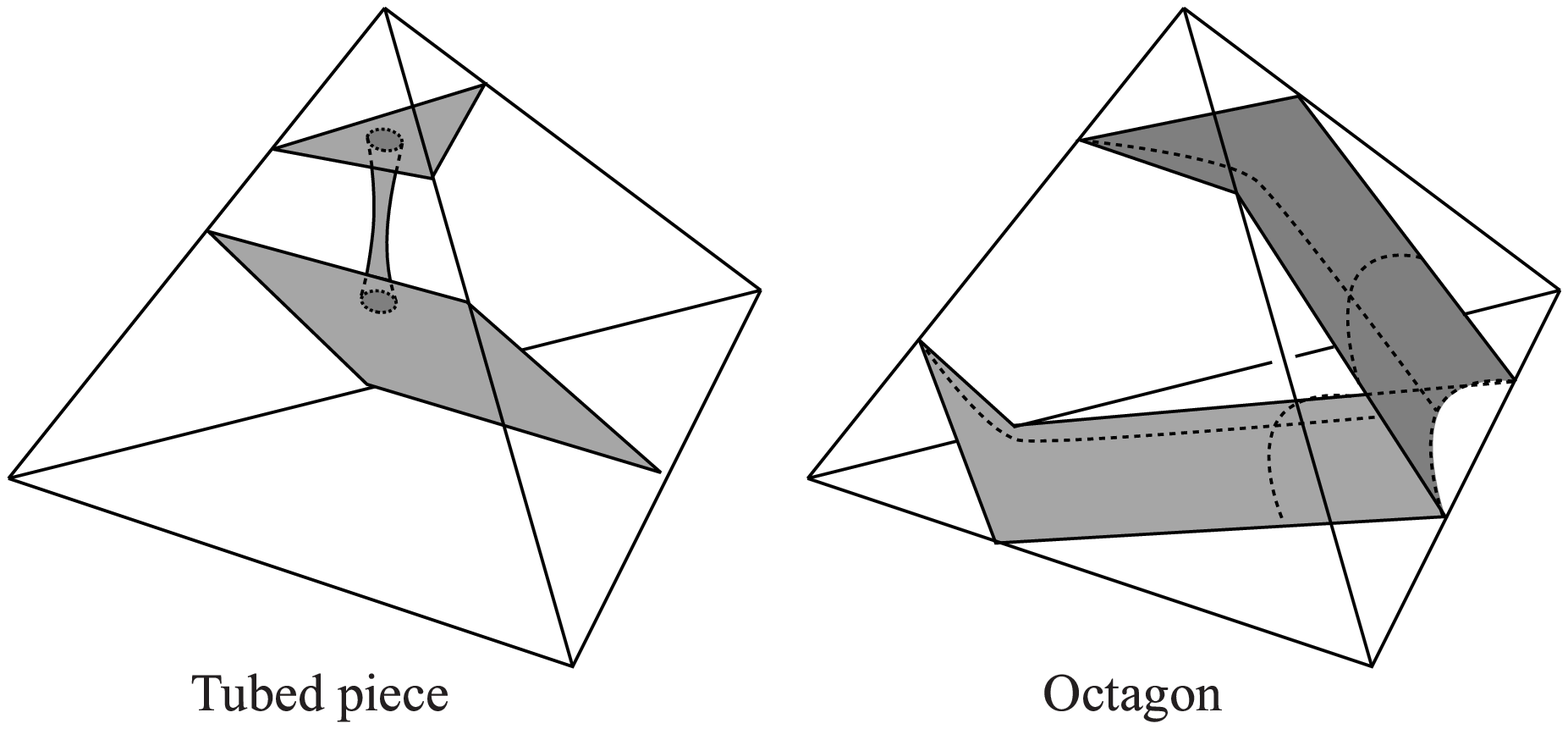}
}
\vskip 6pt
\centerline{Figure 6: almost normal pieces}

It is a theorem of Rubinstein [12] and Stocking [15] that any strongly
irreducible Heegaard surface in a compact orientable 3-manifold
may be ambient isotoped into almost normal form with
respect to any given triangulation or ideal triangulation.
A variant of this result (Theorem 4.2) will be vital in this paper.

We now examine how normal and almost normal surfaces
interact with a partially flat angle structure.
We follow Matveev [9] and term a surface {\sl 2-normal} if it
is closed and embedded and it intersects each ideal tetrahedron in a collection
of normal discs and octagons.
It will be useful to consider 2-normal surfaces
with non-negative Euler characteristic. Here, we have the
following result.

\noindent {\bf Theorem 2.1.} {\sl Let ${\cal T}$ be a partially flat
angled ideal triangulation of a compact orientable 3-manifold $M$. Then any connected closed 2-normal 
surface in ${\cal T}$ with non-negative Euler characteristic is 
normally parallel to a toral boundary component of $M$.}

The key tool in the proof of this is a quantity known as 
combinatorial area, which is assigned to any 2-normal (or more general)
surface in $M$. The partially flat angle structure assigns an interior
angle in the range $[0, \pi]$ to each edge of each ideal tetrahedron.
The corresponding {\sl exterior angle} is defined to be 
$\pi$ minus the interior angle. The {\sl combinatorial area}
of any normal or almost normal piece is defined to be 
the sum of the exterior angles of the edges it runs over (counted with multiplicity)
minus $2 \pi$ times its Euler characteristic. It is easy to
verify that this is always at least zero. Any triangle
running over edges with interior angles that sum to $\pi$
has zero area. The only other normal or almost normal
piece with zero area is a  so-called {\sl vertical
square} in a flat ideal tetrahedron. This is a square that
intersects the edges with angle $\pi$. (See Figure 7.)

\vfill\eject
\centerline{
\epsfxsize=4in
\epsfbox{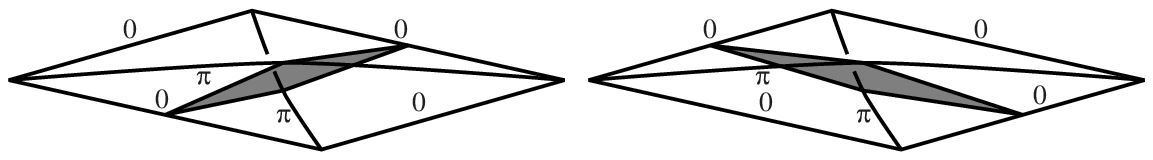}
}
\vskip 6pt
\centerline{Figure 7: vertical squares}

The {\sl combinatorial area} of a 2-normal
surface $F$ is the sum of the combinatorial
areas of its normal and almost normal discs. It is proved
in Proposition 4.3 of [6] that this is equal to $-2\pi \chi(F)$.
Thus, $\chi(F)$ is always non-positive. Suppose, as in the hypothesis of
Theorem 2.1, that $\chi(F)$ is also non-negative. Then, $F$ must be
composed entirely of triangles and vertical squares. We claim
that in fact $F$ consists only of triangles. Hence, if
$F$ is connected, it is normally parallel to a toral boundary component
of $M$, as required.

Vertical squares lie inside flat ideal tetrahedra, each
of which lies in some layered polygon $P$. Focus on the top
ideal tetrahedron of this layered polygon. Its intersection
with the top ideal polygon of $P$ is a square $S$.
One of the diagonals of this square is
an edge in the top ideal polygon. Since no other layered
polygons are attached to this edge, the normal discs
of $F - P$ adjacent to it are all triangles.
Hence, every arc of $S \cap F$
separates off a single vertex of $S$.
Thus, in the intersection of the top ideal tetrahedron
with $F$, there are no vertical squares,
only triangles. Repeating this argument for
each of the ideal tetrahedra of $P$, we deduce that
$F \cap P$ is only triangles, as required. $\square$

It is not known whether any finite-volume hyperbolic 3-manifold admits an 
angled ideal triangulation. But, according to the following existence theorem,
it does always have a partially flat angled ideal triangulation.

\noindent {\bf Theorem 2.2.} {\sl Let $M$ be a compact connected orientable
3-manifold with non-empty boundary. Let $T$ be its toral boundary components.
Then the following are equivalent:
\item{(1)} $M$ is simple and not a 3-ball;
\item{(2)} $M - T$ admits a finite-volume hyperbolic structure with
totally geodesic boundary;
\item{(3)} $M$ admits a partially flat angled ideal triangulation.

\noindent Moreover, if these conditions are satisfied, there is an algorithm
that constructs a partially flat angled ideal triangulation, starting with
any triangulation of $M$.}

\noindent {\sl Proof.} (1) $\Rightarrow$ (2): This is a well known result of
Thurston. The proof goes as follows. Let $DM$ be the result of doubling $M$
along $\partial M - T$, and let $DT$ be the two copies of $T$ in $DM$.
Then $DM$ is a compact orientable simple Haken 3-manifold with (possibly empty)
toral boundary. So, by Thurston's geometrisation theorem [10], $DM - DT$ admits a complete
finite-volume hyperbolic structure. There is an involution of $DM$ that interchanges
its two halves. By Mostow's rigidity theorem, this is homotopic to an isometry.
By a result of Tollefson [16], the involution and the isometry are equivariantly
isotopic. The fixed-point set of this isometry is therefore a totally
geodesic copy of $\partial M - T$ in $DM$. This divides $DM - DT$ into two
copies of $M - T$, each of which inherits a finite-volume hyperbolic
structure with totally geodesic boundary, as required.

(2) $\Rightarrow$ (3): It is a theorem of Epstein and Penner [1]
that, when $\partial M = T$, the interior of $M$ is obtained from a finite collection of hyperbolic
ideal polyhedra, by gluing their faces isometrically in pairs.
When $\partial M$ strictly contains $T$, there is a version of this
theorem, due to Kojima [5]. Instead of hyperbolic ideal polyhedra,
one uses truncated hyperbolic hyperideal polyhedra. Recall that
these are defined as follows. Use the projective model for hyperbolic
space ${\Bbb H}^3$, which is the open ball in projective 3-space. A {\sl  hyperideal
polyhedron} is the intersection of this open ball with a polyhedron $P$,
such that every vertex of $P$ lies outside of ${\Bbb H}^3$, but
where no edge of $P$ lies completely outside of ${\Bbb H}^3$.
Thus, some vertices of $P$ lie on the sphere at infinity of ${\Bbb H}^3$
(these are the {\sl ideal} vertices), and some lie outside
the sphere at infinity (these are the {\sl hyperideal} vertices).
Each hyperideal vertex of $P$ is at the apex of
a cone tangent to the sphere at infinity of ${\Bbb H}^3$. The intersection of this cone
with the sphere at infinity is a circle, which bounds a totally geodesic
plane in ${\Bbb H}^3$. If one truncates the polyhedron along each
of these planes and removes any vertices on the sphere at infinity, 
the result is a {\sl truncated hyperbolic hyperideal
polyhedron}. We permit all the vertices of $P$ to lie on the sphere at infinity,
and so a hyperbolic ideal polyhedron is a special case of
a hyperbolic hyperideal polyhedron and a special case of a truncated
hyperbolic hyperideal polyhedron.

Pick a vertex of each polyhedron $P$ as above, that arises in the
decomposition of $M - T$ into truncated hyperbolic hyperideal polyhedra.
We call this the {\sl coning vertex} of $P$.
The polyhedron $P$ is therefore a cone on this vertex, the base of the cone being
those faces that do not contain the vertex. If we subdivide
each of these faces into triangles, then coning these
off at the coning vertex induces a decomposition of the
hyperideal polyhedron into hyperideal tetrahedra. (See Figure 8.)
Each hyperideal tetrahedron
inherits a set of non-zero interior angles, satisfying condition
(i) in the definition of a partially flat angled ideal triangulation.

\vskip 18pt
\centerline{
\epsfxsize=2.8in
\epsfbox{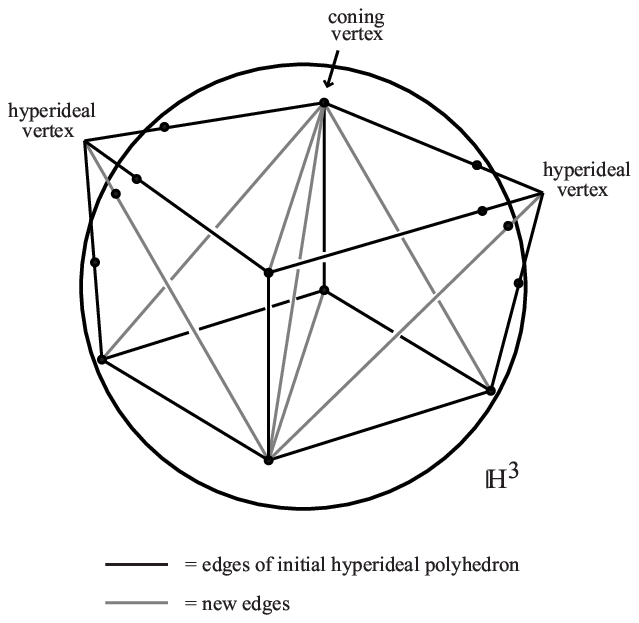}
}
\vskip 6pt
\centerline{Figure 8.}

This decomposition of the hyperideal polyhedron 
yields a decomposition of each its faces into topological ideal
triangles. These topological ideal triangles may not be hyperbolic ideal triangles
because some of their vertices may lie outside the sphere at infinity.
When two faces of the decomposition are glued isometrically,
their topological ideal triangulations may not agree. However,
these two ideal triangulations differ by a finite sequence
of elementary moves, which we may assume leaves no edge in the interior
of the faces untouched. Insert the corresponding layered
polygon between the two faces, interpolating between their
ideal triangulations. (See Figure 9.) Thus, we obtain
an ideal triangulation of the 3-manifold with an angle assignment to each edge
of each ideal tetrahedron. It is clear that
the conditions (i), (ii), (iii) and (iv) in the definition of
a partially flat angled ideal triangulation are satisfied.

(3) $\Rightarrow$ (1): This is essentially contained in Corollary 4.6 in [6].
We sketch the proof now. Suppose that $M$ admits a partially flat
angled ideal triangulation. If $M$ is reducible, then it contains a normal
2-sphere, contrary to Theorem 2.1. If $M$ contains a properly embedded incompressible
torus, then this can be ambient isotoped into normal form, and hence
is boundary parallel by Theorem 2.1. In order to deal with properly embedded
discs and incompressible annuli in $M$, we need to introduce a definition of normal surfaces
that intersect $\partial M$ and to prove a version of Theorem 2.1 for these.
We will not give the full details here, but refer the reader instead to Proposition 4.5 in [6]. 
Thus, $M$ is simple. Also, $M$ cannot be a 3-ball, for one could then find a normal 2-sphere 
parallel to $\partial M$, contradicting Theorem 2.1.

\vskip 18pt
\centerline{
\epsfxsize=4in
\epsfbox{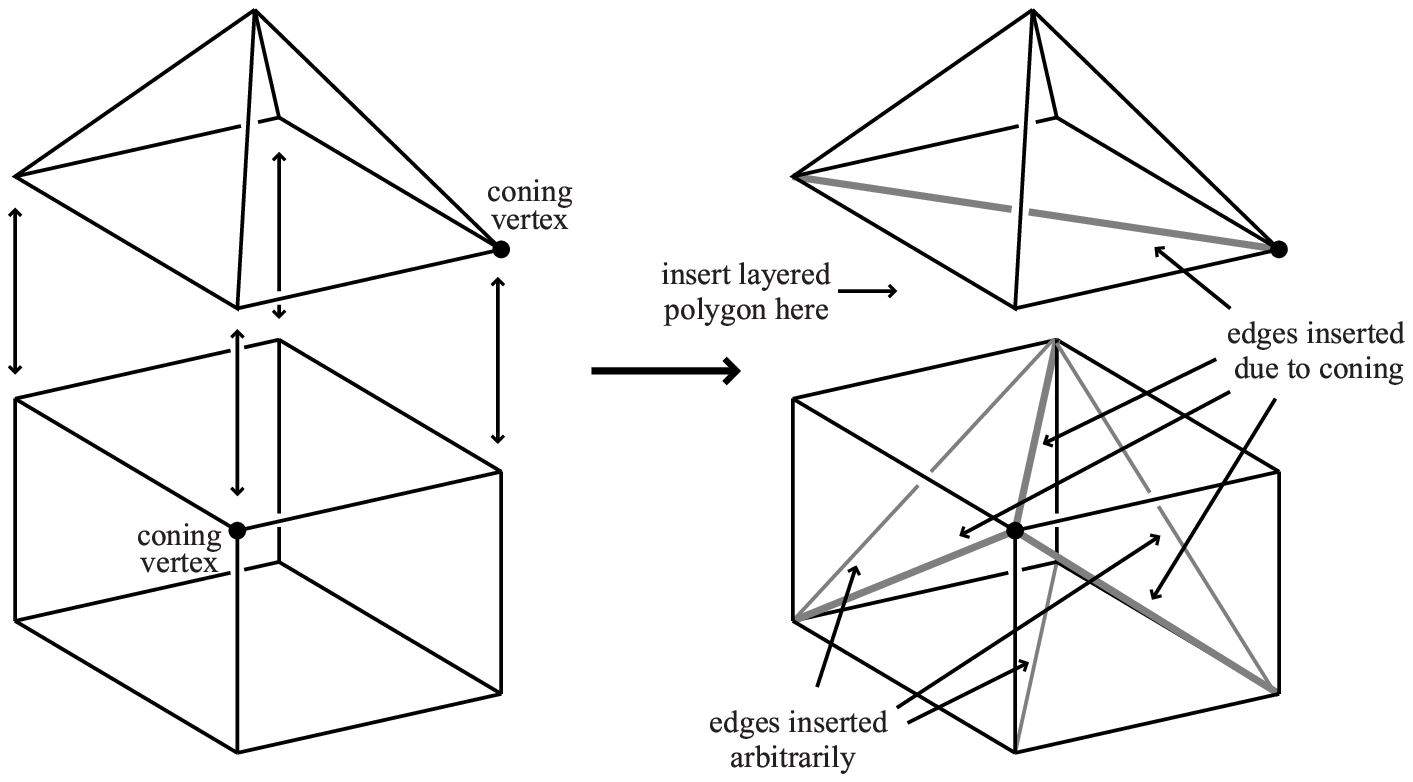}
}
\vskip 6pt
\centerline{Figure 9.}

Let us now suppose that $M$ has a partially flat angled ideal
triangulation. We need to give an algorithm to find one.
Starting with a triangulation of the manifold, there
is a simple algorithm that constructs an ideal triangulation
(see Theorem 1.1.13 of [9]). Given any ideal triangulation, there is an
algorithm that determines whether it admits a 
partially flat angle structure, since this
is just a linear programming problem.
Any two ideal triangulations of a compact orientable
3-manifold differ by a 
sequence of 2-3 and 3-2 moves (see Figure 10), by a result of 
Matveev (Theorem 1.2.5 of [9]).

\vskip 18pt
\centerline{
\epsfxsize=2.3in
\epsfbox{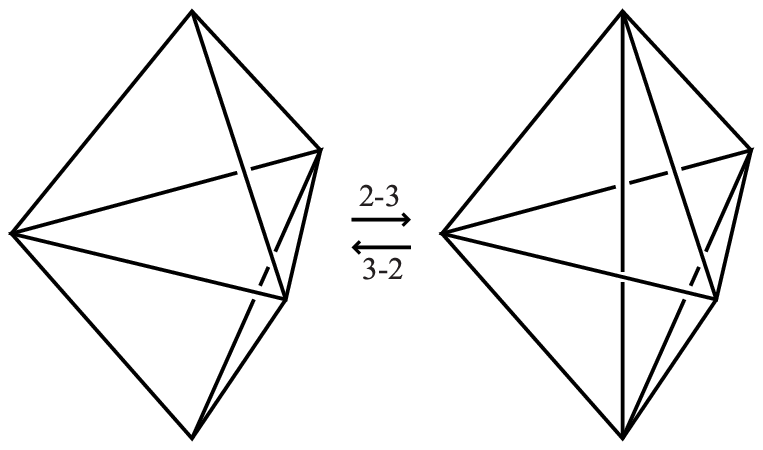}
}
\vskip 6pt
\centerline{Figure 10: $2-3$ and $3-2$ moves}

Thus, the algorithm to construct the partially flat
angled ideal triangulation proceeds as follows. One 
checks whether the initial ideal triangulation admits
a partially flat angle structure. If it does,
we are done and we stop. If not, then one applies
all possible 2-3 and 3-2 moves to the ideal triangulation,
giving a new collection of ideal triangulations.
One checks each of these for partially flat angle
structures. Continuing in this fashion, a partially
flat angled ideal triangulation is eventually constructed.
$\square$

\vskip 18pt
\centerline{\caps 3. Generalised Heegaard splittings}
\vskip 6pt

It is technically convenient, when dealing with Heegaard
surfaces, to focus on those that are strongly irreducible.
The piece of machinery that allows one to make
this reduction is known as untelescoping,
which yields a generalised Heegaard splitting
for the manifold. We now briefly describe these concepts.

Recall that a {\sl compression body} $C$ is a connected
orientable 3-manifold that either is a handlebody or is
obtained from $S \times [0,1]$ by attaching 1-handles
to $S \times \{ 1 \}$, where $S$ is a closed orientable,
possibly disconnected, surface. The copy of $S \times \{0 \}$
in $C$ is termed the {\sl negative boundary} and is
denoted $\partial_- C$. The negative boundary is defined
to be empty when $C$ is a handlebody. The remainder
of $\partial C$ is the {\sl positive boundary} and is
denoted $\partial_+ C$. A {\sl handle structure} on
$C$ is an expression of $C$ as either $\partial_- C \times I$
with 1-handles attached, or as a 3-ball with 1-handles attached.
Note that, in general, a compression body has many different
handle structures.

A {\sl generalised Heegaard splitting} of a compact orientable 3-manifold
$M$ is a decomposition of the manifold along closed orientable disjoint
properly embedded separating surfaces into manifolds
$C_1, \dots, C_m$, each of which is a disjoint
union of compression bodies, such that
$\partial_- C_{2i} \cap {\rm int}(M) = \partial_- C_{2i+1} \cap {\rm int}(M)$
and $\partial_+ C_{2i} = \partial_+ C_{2i-1}$
for each relevant integer $i$.
Let $F_i$ be the surface $C_i \cap C_{i+1}$.
This is known as an {\sl even} or {\sl odd} surface
depending on the parity of $i$.
We view the even surfaces as dividing $M$ into a
collection of 3-manifolds, and the odd surfaces as
forming Heegaard splittings for these manifolds.
(See Figure 11.)

\vskip 18pt
\centerline{
\epsfxsize=1in
\epsfbox{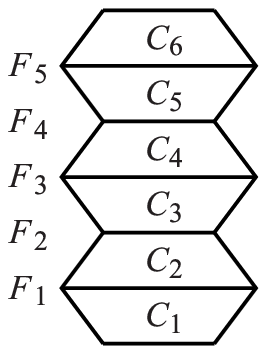}
}
\vskip 6pt
\centerline{Figure 11: a generalised Heegaard splitting}

There is a method for constructing a Heegaard splitting
for a 3-manifold, starting with a generalised Heegaard
splitting $\{ C_1, \dots, C_m \}$, known as {\sl amalgamation} [14]. This procedure
is a sequence of modifications, each of which we term
a {\sl partial amalgamation}. Each partial amalgamation
is based around one of the even surfaces, $F_2$, say.
Either side of this even surface, there are two collections
of compression bodies $C_2$ and $C_3$. Pick a handle
structure on each of these compression bodies that is
not a handlebody. Thus, we view each such compression body as obtained from
$F_2' \times I$, where $F_2' \times \{ 0 \}$ is the relevant components
of $F_2$, by attaching
a collection of 1-handles to $F_2' \times \{ 1 \}$. We extend each of these
1-handles vertically through $F_2' \times I$, so that
they are attached to $F_2$. We may ensure
that the attaching discs of these 1-handles are all disjoint.
Let $F_1'$ be the surface obtained from $F_2$ by attaching these tubes.
It separates $C_1 \cup C_2 \cup C_3 \cup C_4$ into two collections of
compression bodies $C_1'$ and $C_2'$, where $C_1'$ is a copy of $C_1$ with 1-handles
attached, and $C_2'$ is a copy of $C_4$ with 1-handles attached.
We therefore end with a new generalised Heegaard
splitting $\{ C_1', C'_2, C_5, \dots, C_m \}$ for
$M$, which is obtained from the previous one by partial
amalgamation. (See Figure 12.) When this procedure is performed as many times
as possible, the result is a Heegaard splitting, which is
an {\sl amalgamation} of the original generalised Heegaard
splitting.

\vskip 18pt
\centerline{
\epsfxsize=4in
\epsfbox{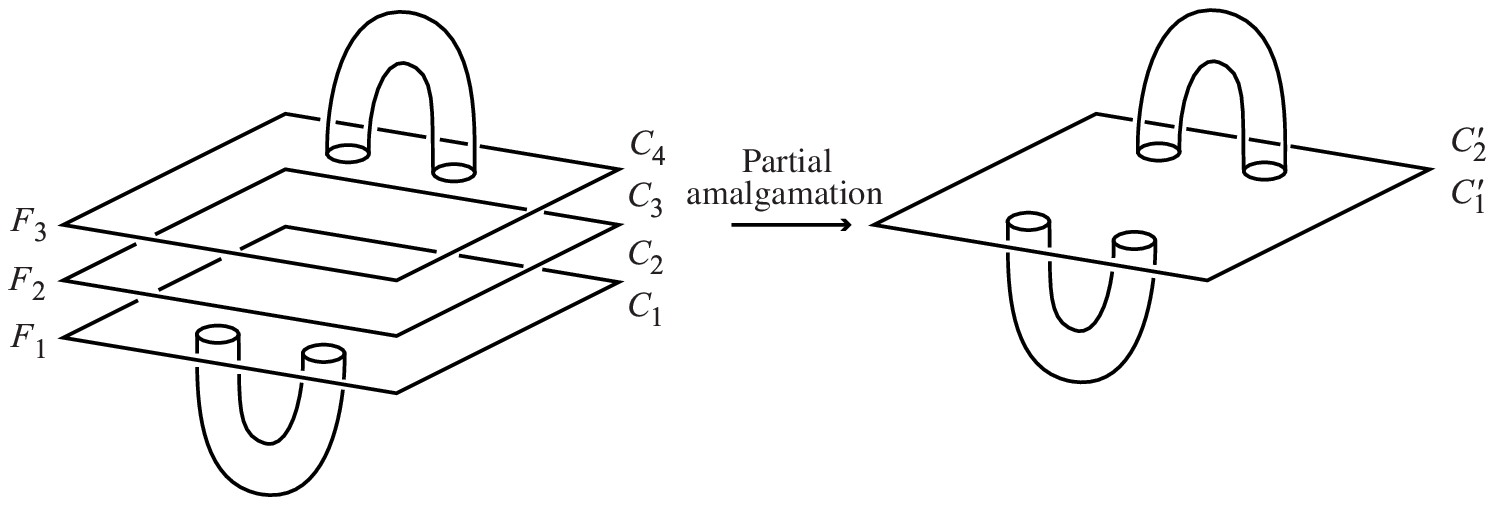}
}
\vskip 6pt
\centerline{Figure 12.}

Choices were made when forming the Heegaard splitting
for $M$: we picked handle structures on $C_2$ and $C_3$,
and we picked an order on the even surfaces in
which to perform the partial amalgamations. It is in fact the
case that the resulting Heegaard splitting of $M$ is
independent of these choices. This important result
does not appear to be present in the literature, and so
we provide a proof.

\noindent {\bf Proposition 3.1.} {\sl If one amalgamates a
generalised Heegaard splitting, the resulting Heegaard
splitting is well-defined up to ambient isotopy. In particular,
it is independent of the order of partial amalgamations and the choice
of handle structures on the compression bodies.}

Let us first examine what happens when we change the order
of the partial amalgamations. Each partial amalgamation is
based around an even surface. So, consider two such even
surfaces, and the associated partial amalgamations. We must
show that if one swaps the order of these partial amalgamations,
the resulting generalised Heegaard splitting is unchanged
up to ambient isotopy.
This is clear if the indexing integers of the even surfaces differ by more
than 2, because in this case none of the compression bodies
involved in the different partial amalgamations intersect. Thus, we focus
on the case where the indexing integers of the even surfaces
differ by 2: say that they are $F_2$
and $F_4$. Now, we may view the former partial amalgamation
procedure as the removal of $F_1$ and $F_2$, and the addition
of handles onto $F_3$. Similarly, the latter partial amalgamation
can be viewed as the removal of $F_4$ and $F_5$, together with
addition of handles onto $F_3$. So, whatever the order of
the two partial amalgamations, the resulting odd surface is the
same: it is $F_3$ with handles attached to both sides.

Let us now consider what happens when we vary the
handle structure on one of the compression bodies $C$ that
is not a handlebody. This handle structure
is determined by the co-cores of the 1-handles, which
form a collection $D$ of disjoint compression discs for
$\partial_+ C$. This collection is {\sl complete}, in
the sense that the result of compressing $\partial_+ C$
along $D$ is a copy of $\partial_- C$. There is clearly
a one-one correspondence between handle structures on
$C$ (up to ambient isotopy) and complete collections 
of compression discs for $\partial_+ C$ (up to ambient
isotopy). Thus, we are led
to the question of how two different complete collections
of compression discs are related. The following answer is
well known (see Corollary 1.6 of [4] for example).

\noindent {\bf Lemma 3.2.} {\sl Any two complete collections
of compression discs for a compression body differ by
a finite sequence of band moves.}

The definition of a {\sl band move} is as follows. Let $D_1$
and $D_2$ be distinct discs in a complete collection $D$. Let $\alpha$ be an arc
in $\partial_+ C$ with interior disjoint from $D$ and
with one endpoint in $D_1$ and the other endpoint in $D_2$.
Let $N$ be a regular neighbourhood
of $D_1 \cup \alpha \cup D_2$. Then,
${\rm cl}(\partial N - \partial_+ C)$ consists of three
compression discs for $\partial_+ C$, one parallel to $D_1$,
one parallel to $D_2$, and a third which we denote by $D_1'$.
Then, $D \cup D_1' - D_1$ is a new complete collection
of compression discs, obtained from $D$ by a band move.

When the compression body $C$ is embedded within a 3-manifold $M$,
as in the current situation, we can realise these band moves by 
{\sl handle slides} as follows. We view $C$ as $\partial_- C \times [0,1]$,
with 1-handles attached to $\partial_-C \times \{ 1 \}$. Let $D = 
\{D_1, D_2, \dots, D_n \}$ be the cocores of the 1-handles, and let $D'
= \{ D'_1, D_2, \dots, D_n \}$ be obtained from $D$ by a band move
along $\alpha$, as above. Now isotope the
1-handle corresponding to $D_2$, by sliding its attaching disc incident
to $\alpha$ along $\alpha$, and then over $D_1$. The new compression body
$C'$ is clearly ambient isotopic to $C$, but now $D'$ is ambient isotopic in $C'$
to the cocores of its 1-handles. See Figure 13.

\vskip 6pt
\centerline{
\epsfxsize=4.4in
\epsfbox{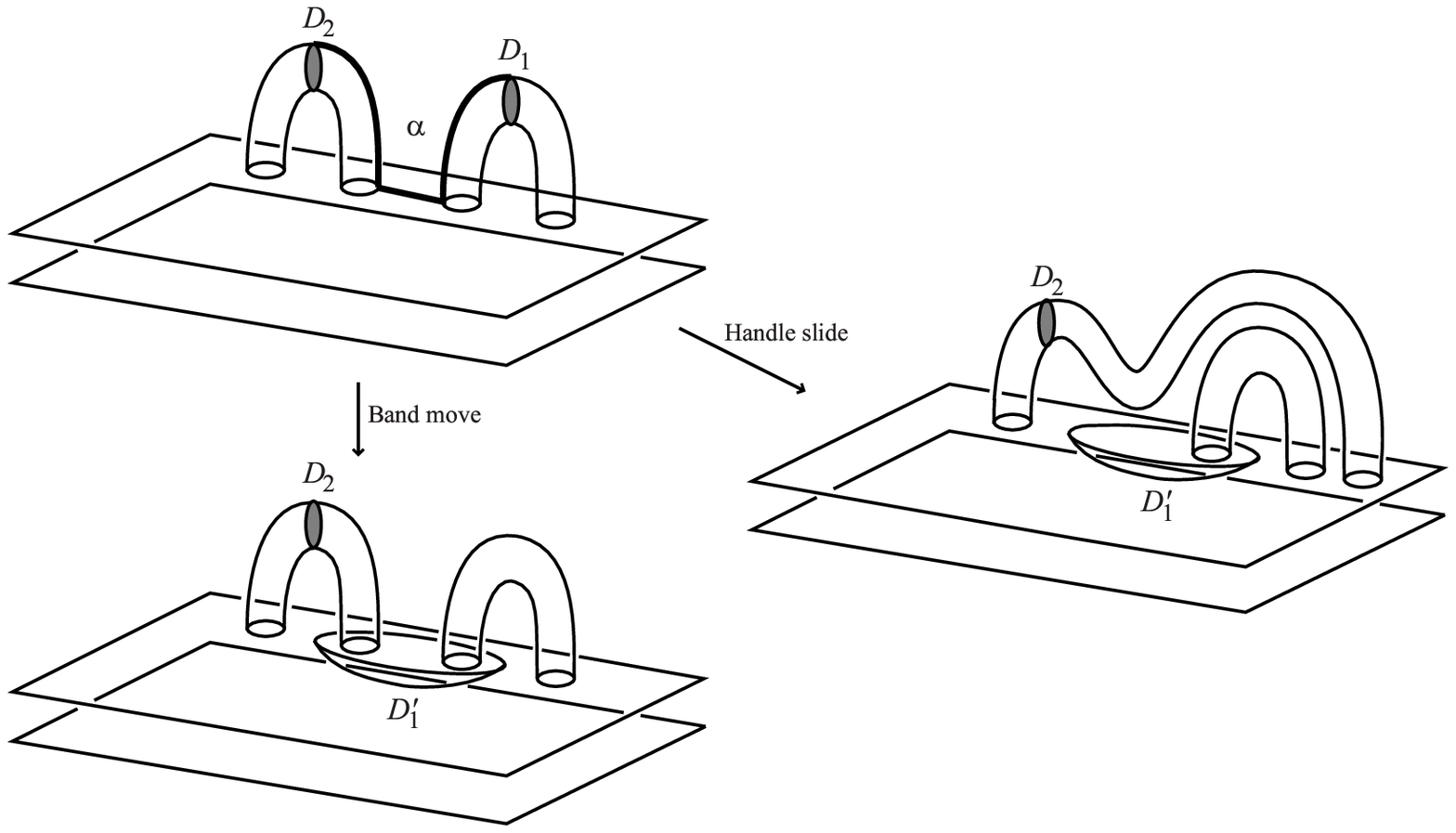}
}
\vskip 6pt
\centerline{Figure 13.}

Let us now consider two different ways of performing a partial
amalgamation upon  a generalised Heegaard splitting $\{ C_1, \dots, C_m \}$.
Let us suppose, for the sake of being definite, that these
partial amalgamations are centred on the surface
$\partial_- C_2 = \partial_- C_3$. Let us also suppose that
these two partial amalgamations are specified by the same
handle structures on $C_2$, but where the handle structures
on $C_3$ differ by a handle slide.
We view $C_3$ as $\partial_- C_3 \times [0,1]$, with 1-handles
attached in two different ways. The second set of
handles is obtained from the first set, by sliding
one of the attaching discs of one of the handles
along an arc $\alpha'$ in $\partial_- C_3 \times \{ 1 \}$ up to another handle and
over that handle. When we perform the first partial amalgamation,
the new odd surface $F_1'$ is obtained from $\partial_- C_3$ by attaching
handles onto both sides. By using the product structure on
$\partial_- C_3 \times [0,1]$, we may project the
arc $\alpha'$ to an arc in $\partial_- C_3$. By slightly isotoping
this arc if necessary, we may assume that its interior avoids the
attaching discs of all the handles. This
arc joins two handles of $F_1'$. We may therefore slide one
of these handles along this arc, and over the other handle.
The resulting surface $F_1''$ is exactly that obtained by the
second partial amalgamation. Thus, $F_1'$ and $F_1''$ are related by a handle
slide and are therefore ambient isotopic. (See Figure 14.)

We have therefore shown that the choices made in creating the amalgamated
Heegaard surface do not affect its ambient isotopy class. This proves
Proposition 3.1. $\square$

\vskip 18pt
\centerline{
\epsfxsize=4.5in
\epsfbox{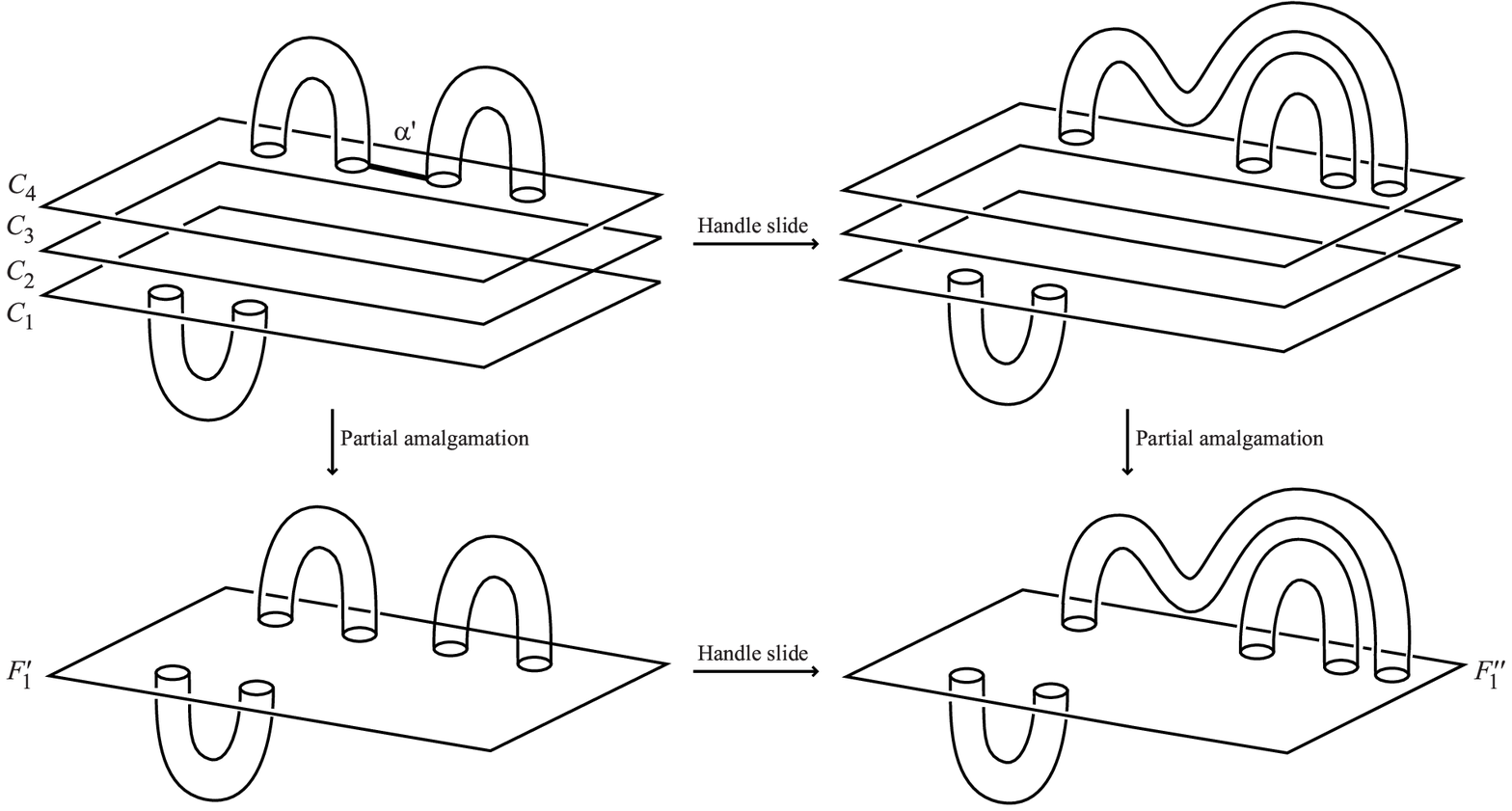}
}
\vskip 6pt
\centerline{Figure 14.}

We will be constructing the Heegaard splittings 
required by Theorem 1.1 by first constructing generalised
Heegaard splittings. Thus, we need to know that the process
of amalgamation can be achieved algorithmically.

\noindent {\bf Proposition 3.3.} {\sl Let $F$ be a Heegaard
surface for $M$ that is obtained from a generalised Heegaard
splitting $\{ C_1, \dots, C_m \}$ by amalgamation. Suppose that
$\{ C_1, \dots, C_m \}$ is given as a subcomplex of a triangulation of $M$. Then,
there is an algorithm that constructs $F$ in $M$.}

\noindent {\sl Proof.} Let $\{ C'_1, C'_2, C_5, \dots, C_m \}$
be obtained from $\{ C_1, \dots, C_m \}$ by a partial amalgamation.
It clearly suffices to construct
$C'_1$ and $C'_2$ from $C_1$, $C_2$, $C_3$ and $C_4$.
Let us focus on a component of $C_3$, say, that is not a handlebody. A complete
set of compression discs for this compression body is
constructible (see Theorem 4.1.14 of [9] or Algorithm 9.3 of [3]). Cutting along this collection, we obtain
a copy of $F_2' \times [0,1]$, where $F_2'$ is the relevant components of
$\partial_- C_3$. In $F_2' \times \{ 1 \}$, we
have two copies of each compression disc, giving a
collection $D'$ of disjoint discs. We may construct
$D' \times [0,1]$ in $F_2' \times [0,1]$ as follows. In $F_2' \times [0,1]$,
we may construct a vertical annulus $A$, using Theorem 6.4.10 of [9] or Algorithm 9.7 of [3]. 
By performing an ambient isotopy on $A$ supported in a small
neighbourhood of $F_2' \times \{ 1 \}$, we may ensure that
$A$ intersects each component of $D'$ in a non-empty collection
of arcs. For each disc $D''$ in $D'$, we may construct a properly embedded arc in $A$
running from a component of $D'' \cap A$ to $F_2' \times \{ 0 \}$. We may
arrange that these arcs are pairwise disjoint. Let $R$ be
their union. Then, a regular neighbourhood of $R \cup D'$
in $F_2' \times [0,1]$ is the
required copy of $D' \times [0,1]$. If we now reverse the cutting
procedure which gave $F_2' \times [0,1]$ from the component of $C_3$,
the components of $D' \times [0,1]$ glue up in pairs to form
a collection of 1-handles attached to $\partial_- C_3$. Perform this procedure 
for each compression body component of $C_2$ and $C_3$,
and then isotope if necessary, so that the attaching
discs of the 1-handles in $\partial_- C_2 = \partial_- C_3$
are disjoint. The new surface $\partial_- C_1'
= \partial_- C'_2$ is obtained from $\partial_- C_2$ by attaching these tubes.
$\square$

The following result was proved by Scharlemann and Thompson [13]. It
describes a process known as {\sl untelescoping}.

\noindent {\bf Theorem 3.4.} {\sl Let $M$ be a compact orientable
3-manifold, and let $F$ be an irreducible Heegaard splitting.
Then there is a generalised Heegaard splitting $\{ C_1, \dots, C_m \}$ for $M$,
such that 
\item{(i)} the even surfaces are incompressible and have
no 2-sphere components;
\item{(ii)} the odd surfaces are strongly irreducible;
\item{(iii)} no $C_i$ is homeomorphic to $\partial_- C_i \times I$ 
(although some components of $C_i$ may be products);
\item{(iv)} $F$ is obtained from this generalised Heegaard splitting by
amalgamation.

\noindent Suppose, in addition, that the Heegaard genus of $M$ is
more than $1$. Then, we may also arrange that no odd surface is
composed entirely of tori.}

We now wish to estimate the genus of the odd and even surfaces in this
generalised Heegaard splitting. Let us suppose that
the Heegaard genus of $M$ is more than $1$.
Now, it is trivial to check that the quantity
$$\sum_{i=1}^m {\chi(\partial_- C_i)  - \chi(\partial_+ C_i) \over 2}$$
is unchanged under partial amalgamation. Hence, it equals $-\chi(F)$.
Each term in the sum is a positive integer, by (iii) and the fact that no $C_i$
is a collection of solid tori, and no component of any $C_i$ is a 3-ball. Thus, we obtain
the inequality $m \leq -\chi(F)$. Since $F$ is obtained from
the splitting by amalgamation, it can be viewed as obtained from
any given even or odd surface by adding tubes. Thus, the genus
of each even or odd surface is at most $g(F)$, the genus of $F$. The number
of even and odd surfaces is $m-1 \leq 2g(F) - 3$. So, the genus
of the union of the odd and even surfaces is at most
$g(F)(2g(F)-3)$. (It is possible to improve this estimate slightly,
but all that is needed here is a computable upper bound on the
genus of the union of the odd and even surfaces in terms of $g(F)$.)

\vfill\eject
\centerline{\caps 4. Almost normal surfaces}
\vskip 6pt

This paper relies heavily on the following important theorem of
Rubinstein [12] and Stocking [15].

\noindent {\bf Theorem 4.1.} {\sl Let $M$ be a compact orientable
irreducible 3-manifold, with a given triangulation. 
Let $F$ be a strongly irreducible Heegaard surface for $M$.
Then there is an ambient isotopy taking $F$ into almost
normal form.}

In this paper, we need the following slight extension of
this result, which deals also with ideal triangulations
and with generalised Heegaard splittings.

\noindent {\bf Theorem 4.2.} {\sl Let $M$ be a compact orientable
irreducible 3-manifold, with a given triangulation or ideal
triangulation. Let $\{ C_1, \dots, C_m \}$ be a generalised Heegaard splitting
for $M$. Suppose that the even surfaces are incompressible and
have no 2-sphere components and the odd surfaces are strongly
irreducible. Then there is an ambient isotopy that makes each 
even surface normal and each component of the odd surfaces almost normal.}

The proof follows the argument of Stocking in [15]
almost word-for-word. We refer the reader to [15] for more
details.

For our purposes here, the main usefulness of normal and almost normal
surfaces is that they are constructible.

\noindent {\bf Theorem 4.3.} {\sl Let ${\cal T}$ be a partially flat angled
ideal triangulation of a compact orientable 3-manifold $M$. Then, for any integer $n$, ${\cal T}$ contains
only finitely many closed orientable properly embedded surfaces $F$ with 
${\rm genus}(F) \leq n$, and where each component of $F$ is 
either normal or almost normal. Moreover, there is
an algorithm to construct each of these surfaces.}

Let $F$ be a closed orientable properly embedded surface, each component of
which is normal or almost normal. Let $\overline F$ be
obtained from $F$ by compressing any tubed pieces.
Thus, $\overline F$ is 2-normal, and ${\rm genus}(\overline F)
\leq {\rm genus}(F) \leq n$. So, it clearly suffices to
construct a finite list of possibilities for $\overline F$.
For we may then reconstruct $F$ by reattaching tubes running
parallel to the edges of ${\cal T}$. Note that,
according to Theorem 2.1, ${\cal T}$ contains no 2-normal 2-spheres. 
Hence, each of the compressions we performed on $F$ was essential.

Now, $\overline F$ may be specified
by a vector, each co-ordinate of which is a non-negative
integer, as follows. One associates with each ideal tetrahedron
10 co-ordinates. Each co-ordinate corresponds to a type
of normal or almost normal disc in that tetrahedron:
4 triangle types, 3 square types and 3 octagon types. 
The vector corresponding
to $\overline F$ simply counts the number of copies in $\overline F$ of each
normal and almost normal disc in each ideal tetrahedron.
The fact that these discs patch together to form a
closed surface forces this vector to satisfy certain
linear equations. There are three equations for each
face of the ideal triangulation, corresponding to the three 
types of properly embedded arc in that face. These are
known as the {\sl matching equations}. An embedded surface
cannot contain different square or octagon types
in any given ideal tetrahedron. This again forces constraints
on the vector of $\overline F$. For normal surfaces, these
are known as the {\sl quadrilateral conditions}. In
our situation, we will term them the {\sl square/octagon
conditions}. There is a one-one correspondence between closed
properly embedded 2-normal surfaces and non-negative integer solutions
to the matching equations that satisfy the square/octagon
conditions.

Crucial is the concept of {\sl normal sum}.
Suppose that the vector corresponding to $\overline F$
can be written as a sum of vectors, each of which
has non-negative integer co-ordinates and satisfies
the matching equations. Then these
vectors also satisfy the square/octagon conditions
and so correspond to 2-normal surfaces $F_1$
and $F_2$. We write $\overline F = F_1 + F_2$. It is easy 
to check that $\chi(\overline F) = \chi(F_1) + \chi(F_2)$. 
When $\overline F$ cannot be written as a sum
of non-empty 2-normal surfaces, $\overline F$
is said to be {\sl fundamental}. Crucial to our
algorithms is the following fact (see Theorem 3.2.8 of [9]).

\noindent {\bf Lemma 4.4.} {\sl There is a finite computable list of
fundamental surfaces, such that any 2-normal surface may
be written as a sum of these fundamental surfaces.}

Denote these fundamental surfaces by $F_1, \dots, F_m$.
Suppose that $F_1, \dots, F_r$ are normally parallel to
toral boundary components of $M$ and that the rest are not.
Thus, Lemma 4.4 states that any 2-normal surface $\overline F$
can be written as $\sum_{i=1}^m n_i F_i$, for non-negative 
integers $n_i$. Consider $\sum_{i={r+1}}^m n_i F_i$, which
is a solution to the matching equations satisfying
the square/octagon conditions. It therefore corresponds to
a 2-normal surface $F'$. According to Theorem 2.1, $\chi(F_i)$
is strictly negative for each $i > r$. Hence, it is at most $-1$,
and we obtain the inequalities 
$$\sum_{i=r+1}^m n_i \leq 
-\sum_{i=r+1}^m n_i \chi(F_i) = 
- \sum_{i=1}^m n_i \chi(F_i) =
- \chi(\overline F)  = 2g(\overline F) - 2|\overline F| < 2g(\overline F)
\leq 2n.$$ 
Thus, there is a finite list of
possibilities for $F'$ and they are all constructible.

The surface $\sum_{i=1}^r n_i F_i$
is a collection of copies of the toral boundary components, which
we may realise as disjoint from $F'$. Thus, the union of these
surfaces and $F'$ is a solution to the matching equations
with the same vector as $\overline F$. They are therefore
ambient isotopic. In other words, 
$$\overline F = F' \cup \bigsqcup_{i=1}^r n_i F_i.$$
Hence,
$$g(\overline F) = g(F') + \sum_{i=1}^r n_i.$$
Since, we are assuming that the genus of $\overline F$ is
at most $n$, this provides an upper bound on $\sum_{i=1}^r n_i$.
Thus, there is a finite list of possibilities for $\overline F$
and they are all constructible. The same is then true for $F$.
This completes the proof of Theorem 4.3. $\square$

\vskip 18pt
\centerline{\caps 5. The algorithms}
\vskip 6pt

We now have all the ingredients to describe the algorithms in Theorem 1.1
and to prove that they work. Note that the first algorithm, which computes
the Heegaard genus of $M$, can be constructed from the second algorithm,
which finds all Heegaard surfaces with genus at most a given integer $n$.
This is done as follows.
One first sets $n$ to be $2$ (the smallest possible Heegaard genus
for $M$) and one searches for Heegaard surfaces with genus at most $n$.
If there is one, the Heegaard genus is $2$. If there is not,
set $n$ to be $3$, and repeat.  The first time the algorithm
finds a Heegaard surface, it necessarily has minimal genus,
and the algorithm stops.

Therefore, let us fix a non-negative integer $n$. We will
describe the algorithm to find all Heegaard surfaces in
$M$ with genus at most $n$.

We may restrict attention to irreducible Heegaard surfaces.
For if a Heegaard surface in $M$ is reducible, it is stabilised,
and is therefore obtained from an irreducible Heegaard surface
of smaller genus by stabilising a number of times.

\noindent {\bf Step 1.} Find a partially flat angled ideal
triangulation for $M$.

The algorithm to achieve this is described in the proof of Theorem 2.2.
The algorithm finds not just the required ideal triangulation with an explicit
partially flat angle structure,
but also provides a method of constructing it from
the initial given triangulation.

\noindent {\bf Step 2.} Finding candidates for generalised Heegaard
splittings.

According to Theorem 3.4, given any irreducible Heegaard surface $F$ for $M$,
there is a generalised Heegaard splitting,
from which $F$ is obtained by amalgamation, and in
which each even surface is incompressible
and has no 2-sphere components, and each odd surface is
strongly irreducible. Let $F'$ be the union of the even
and odd surfaces. As observed at the end of section 3,
the conclusions of Theorem 3.4 imply that the genus of $F'$ is at most $g(F)(2g(F) - 3)$.
By Theorem 4.2, we may make each
even surface normal and each component of the odd
surfaces almost normal. According to Theorem 4.3,
there is an algorithm that constructs a finite list
of surfaces in $M$, one of which is $F'$. Step 2 in the
algorithm is to construct this list of surfaces.

\noindent {\bf Step 3.} Determining which are generalised Heegaard splittings.

There is an algorithm to determine whether a properly embedded closed, 
possibly disconnected, surface $F'$ forms a generalised Heegaard splitting.
It proceeds as follows. Cut $M$ along $F'$. There is an algorithm
that determines whether each component of the complement is a compression
body (see Theorem 4.1.14 of [9] or Algorithm 9.3 of [3]).
If this holds, the algorithm then checks all possible
ways of grouping these compression bodies into an ordered collection
$\{ C_1, \dots, C_m \}$ (where each $C_i$ may be disconnected)
such that
$\partial_- C_{2i} \cap {\rm int}(M) = \partial_- C_{2i+1} \cap {\rm int}(M)$
and $\partial_+ C_{2i} = \partial_+ C_{2i-1}$
for each relevant integer $i$.
We apply this algorithm to each surface provided by Step 2, and thereby
create a list of generalised Heegaard splittings.

\noindent {\bf Step 4.} Amalgamation.

Consider one of generalised Heegaard splittings $\{ C_1, \dots, C_m \}$
in our list, and suppose that $F$ is the Heegaard surface obtained
from this by amalgamation. According to Proposition 3.1, this surface
$F$ depends only on the generalised Heegaard splitting, and not on any
choices made during the amalgamation procedure. If one is interested
only in the genus of $F$, then this can be calculated from the
surfaces in the generalised Heegaard splitting via the formula
$$-\chi(F) = \sum_{i=1}^m {\chi(\partial_- C_i)  - \chi(\partial_+ C_i) \over 2}.$$
Thus, if one is interested only in the existence of a Heegaard surface
with genus at most $n$, then this can be determined by applying this
formula to each generalised Heegaard splitting in the list.
The algorithm discards all those Heegaard surfaces with genus more than $n$.
However, if one actually wants to construct all such Heegaard
surfaces, then one must perform each amalgamation algorithmically,
using Proposition 3.3. The result is a finite list of Heegaard surfaces
for $M$ with genus at most $n$. $\square$

\vfill\eject
\centerline{\caps References}
\vskip 6pt

\item{1.} {\caps D. Epstein, R. Penner}, {\sl 
Euclidean decompositions of noncompact hyperbolic manifolds.}
J. Differential Geom. 27 (1988) 67--80. 

\item{2.} {\caps W. Jaco, J. H. Rubinstein}, {\sl
Layered-triangulations of 3-manifolds}, Preprint, 
arxiv:math.GT/0603601

\item{3.} {\caps W. Jaco, J. Tollefson,} {\sl
Algorithms for the complete decomposition of a closed $3$-manifold.}
Illinois J. Math. 39 (1995) 358--406.

\item{4.} {\caps K. Johannson,} {\sl Topology and combinatorics of $3$-manifolds.}
Lecture Notes in Mathematics, 1599. Springer-Verlag (1995).

\item{5.} {\caps S. Kojima}, {\sl Polyhedral decomposition of hyperbolic $3$-manifolds
with totally \break geodesic boundary}, Aspects of low-dimensional manifolds, 93--112,
Adv. Stud. Pure Math. 20, Tokyo (1992).

\item{6.} {\caps M. Lackenby}, {\sl 
Word hyperbolic Dehn surgery}, Invent. Math. 140 (2000) 243--282. 

\item{7.} {\caps T. Li,} {\sl 
Heegaard surfaces and measured laminations, I: the Waldhausen conjecture.}
Invent. Math. 167 (2007) 135--177.

\item{8.} {\caps T. Li,} {\sl
Heegaard surfaces and measured laminations, II: non-Haken 3-manifolds.}
J. Amer. Math. Soc. 19 (2006) 625--657.

\item{9.} {\caps S. Matveev}, {\sl Algorithmic topology and classification of 3-manifolds.}
 Algorithms and Computation in Mathematics, 9. Springer-Verlag, Berlin, 2003. 

\item{10.} {\caps J. Morgan}, {\sl The Smith conjecture.} 
Pure Appl. Math. 112 (1984).

\item{11.} {\caps C. Petronio, J. Weeks}, 
{\sl Partially flat ideal triangulations of cusped hyperbolic 3-manifolds.}
Osaka J. Math. 37 (2000) 453--466. 

\item{12.} {\caps J. H. Rubinstein,} {\sl 
Polyhedral minimal surfaces, Heegaard splittings and decision problems for $3$-dimensional manifolds.} 
Geometric topology (Athens, GA, 1993), 1--20,
AMS/IP Stud. Adv. Math., 2.1, Amer. Math. Soc., Providence, RI, 1997. 

\item{13.} {\caps M. Scharlemann, A. Thompson,} {\sl Thin position for $3$-manifolds.}
Geometric topology (Haifa, 1992), 231--238, Contemp. Math., 164, Amer. Math. Soc., 1994.

\item{14.} {\caps J. Schultens,} {\sl The classification of Heegaard splittings for 
(compact orientable surface)$\,\times\,S\sp 1$,} 
Proc. London Math. Soc. 67 (1993) 425--448. 

\vfill\eject
\item{15.} {\caps M. Stocking}, {\sl
Almost normal surfaces in $3$-manifolds.}
Trans. Amer. Math. Soc.  352  (2000) 171--207.

\item{16.} {\caps J. Tollefson,} {\sl Involutions of sufficiently large $3$-manifolds.}
Topology 20 (1981) 323--352.

\vskip 12pt
\+ Mathematical Institute, University of Oxford, \cr
\+ 24-29 St Giles', Oxford OX1 3LB, United Kingdom. \cr

\end